\documentclass{article}
\usepackage{amsmath}
\usepackage{url}
\usepackage{graphicx,psfrag}
\begin{document}
\title{On Finiteness in the Card Game of War}
\author{Evgeny Lakshtanov and Vera
Roshchina}
\date{}
\maketitle

\begin{abstract}
The game of war is one of the most popular international children's card games. In the
beginning of the game, the pack is split into two parts, then on each move the players reveal their top
cards. The player having the highest card collects both and returns them to the bottom of his hand. The
player left with no cards loses. It is often wrongly assumed that this game is deterministic and the
result is set once the cards have been dealt.However, it is not quite so; as the rules of the game do not
prescribe the order in which the winning player will put his take to the bottom of his hand: own card,
then rival's or vice versa: rival's card, then own. We provide an example of a cycling game with fixed
rules. Assume now that each player can seldom but regularly change the returning order. We have
managed to prove that in this case the mathematical expectation of the length of the game is finite. In
principle it is equivalent to the graph of the game, which has got edges corresponding to all acceptable
transitions, having got the following property: from each initial configuration there is at least one path
to the end of the game.
\end{abstract}

\section{The Game of War}

The Game of War is an internationally popular children's card game with a very basic set of rules.  At the start, the pack is shuffled and divided into two equal parts; then each player reveals their top card.  The player with the highest value card then collects both cards and returns them to the bottom of their hand.  This process continues until one player loses all cards.

Many who play this game in childhood lose patience before the
end as the back and forth nature of card transfer often allows
a player to almost reach the point of victory only to see their
hand reduced dramatically in a matter of moments. In effect
children playing are conducting basic mathematical experiments
and observing the development of a dynamic system. It is often
wrongly assumed that this game is deterministic and the result
is set once the cards have been dealt.  However this is not so;
the rules of the game do not stipulate in which order the
winner of each play returns the cards to the bottom of their
hand –- own card first and then rival's or vice versa.

We shall at first consider a model with an arbitrary (even) number of cards and only one suit.  Hence, a situation where both players present the same value card cannot occur, in the classic game, if this does occur, it would lead to a process called `war' and players would continue to lay cards until one plays a winner, thus claiming all cards played.

\begin{figure}[h]
\centering
\includegraphics[height=100pt, keepaspectratio]{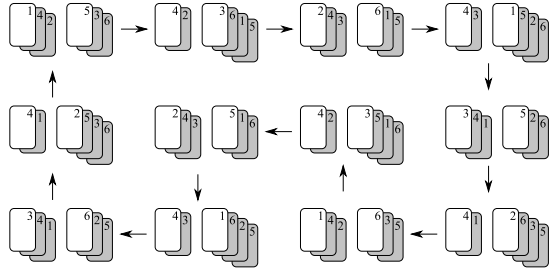}
\caption{ A never ending game, $n=6$.}\label{fig:cycle6}
\end{figure}

\begin{figure}[h] \centering
\includegraphics[height=60pt, keepaspectratio]{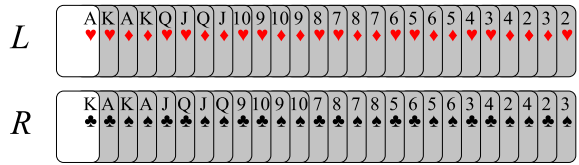}
\caption{ A never ending game, standard pack of 52 cards.}\label{fig:seq}
\end{figure}

So consider a game where the players strictly control the order in which the cards are returned to the bottom of their hand; in this case there is a chance that the game will never finish.  Such a cyclic game is shown in Fig.~\ref{fig:cycle6} where $n = 6$. The card of the `left' player is always returned to the bottom of the winner's hand before the card of the `right' player.
%\newpage
Such a never ending game can also occur with a standard deck if
the above rules are followed.
%\footnote{Note that in the case of
%the standard pack it can happen that on a certain move both
%players show cards of the same value: this situation is called
%``battle''. There exist different versions of rules to resolve
%a battle, but the basic idea is to keep revealing the next pair
%of cards until one of then is higher or one of the players has
%exhausted their hand.}
An example of an initial shuffling for such a game is given on
Fig.~\ref{fig:seq} (throughout the paper we assume that the ace
is the highest value card). Observe that after the first two
moves are made, the last two cards in the hand of player $L$
are $(A\heartsuit, K\clubsuit)$, and in the hand of player
$R$ -- $(K\heartsuit,A\clubsuit)$, then after the next two
moves the last 4 cards are $(A\heartsuit,K\clubsuit,
A\diamondsuit,K\spadesuit)$ for the player $L$ and
$(K\heartsuit,A\clubsuit, K\diamondsuit,A\spadesuit)$ for
player $R$, and so on. It can be easily seen that the order
of the card values is preserved, and after 26 moves the players
will end up with precisely the same distribution of card values
in their hands (although suits would get shuffled). Note that
in this particular case `war' never happens, because the cards
played on each move always have different values.

%%  It should be noted, for the sake of curiosity, that a more complex version of the game can be played, where the highest value card is beaten by the lowest valued, i.e. in the model game with cards numbered from $1$ to $n$, card $n$ beats all other cards except for $1$. Then in the same initial setting either of the players can win. This version of the game is popular in Russia (where 2 wins over ace, but ace beats any card other than 2).  An example of such a game for a model pack of $4$ cards is shown in Fig.~\ref{fig:twooutcomes}.
%%  \begin{figure}[h]
%%  \centering
%%  \includegraphics[height=60pt, keepaspectratio]{Figure7.eps}
%%  \caption{A game with two possible outcomes: $n=4$,  card $4$ wins over all numbers except for $1$}\label{fig:twooutcomes}
%%  \end{figure}

Hence we have established that when rigid rules are used, it might be possible to never finish the game. But what if players use both possible ways of returning cards to their hands, and on each particular move choose such rules at random? In this case, it is not immediately clear whether the players have a nonzero chance to reach the end of the game, however, it is possible to give an answer to this question using a probabilistic model of the game. We show that if the players use both rules to return the cards to their hand, the mathematical expectation of the number of moves in the game is finite, (i.e., there is a zero chance of never finishing the game).

\section{Mathematical model}

\begin{figure}[h] \centering
\includegraphics[height=90pt, keepaspectratio]{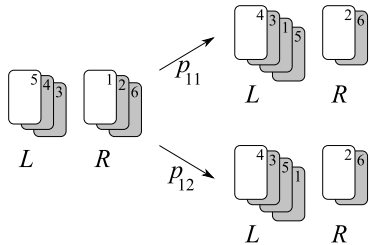}
\caption{Probabilities with which player $L$ chooses the order to place the cards in his hand (model pack of $6$ cards)}\label{fig:prob}
\end{figure}

We now focus on a model case in which there are cards from 1 to
$n$.  We assume that each player possesses certain
peculiarities which means when collecting cards their card will
be placed on the bottom with probability $p^i_1>0$ and second
from bottom with $p^i_2>0$, where the index $i$ identifies the
player ($i$ is either $L$ or $R$).  This is illustrated on
Fig.~\ref{fig:prob} where player $L$ wins, and chooses to place
his card on the bottom with probability $p^{L}_1$, and uses
the other rule with probability $p^{L}_2$.
Let $C$ be the deck of cards, either numbered (from $1$ to $n$, where $n$ is any even number) or standard (with 52 cards, valued from 2 to ace and having 4 suits). We use $L$ to denote the cards in the hand of the `left' player, and $R$ to denote the `right' player's hand. In the model case, $L$ and $R$ are ordered sets of numbers, one of which can be empty; this corresponds to the end of the game.

Each division of the deck into two ordered sets is called a game {\em state}. That is, the game starts with a state with both sets $L$ and $R$ having an equal number of members, and ends with the state in which one of the sets is empty (the {\em final} state). Each game play (when the players show their top cards, compare them and then put the ordered pair to the bottom of the winner's hand) is a {\em transition} from one state to another, because it starts with two ordered sets of cards, and ends with two (different) ordered sets, which correspond to the players' hands.

\begin{figure}[h] \centering
\includegraphics[height=80pt, keepaspectratio]{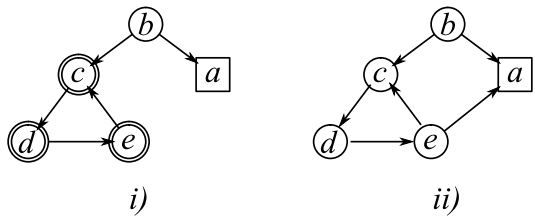}
\caption{Difference between absorbing and non-absorbing graphs: the graph {\it i)} is not absorbing, but the graph {\it ii)} obtained by adding one extra edge between $e$ and $a$ to the graph {\it i)} is absorbing.}\label{fig:graph}
\end{figure}

Such dynamics in which a transition to the next state happens
with fixed probabilities independent on the preceding choices
is called Markov chain.  In the theory of Markov chains, one of
the fundamental facts is the following. Assume that any initial
state is possible.  The mathematical expectation of the number
of moves before the {\em absorption} (reaching the final state)
is finite if and only if the final state can be reached from
any state \cite[Chap.3]{fmc}.\footnote{Let us explain this
fact. Assume that termination is possible from every initial
state. It means that for every vertex $v_i$ there exists $N_i$
such that the conditional probability $P(\cdot | v_i)$ to
terminate the game during next $N_i$ steps is positive. Since
our graph is finite, taking $N=\max\{N_i\}$ and denoting $q>0$ the minimum over i of such probabilities, we see that the
unconditional probability to terminate the game in at most $N$ steps is at least $q>0$. So the probability to
stay in the game $kN$ steps is at most $(1-q)^k$.
From this, the probability to stay in the game $n$ steps converges to zero faster than a
geometric progression. } Usually a Markov
chain is represented as a directed graph in which the vertices
correspond to states, and edges correspond to the transitions.
An edge leaves one vertex and reaches another if and only if
there exists a transition from the former to the latter with a
nonzero probability. It is not difficult to see that in our
case each non-final vertex (or state)  has only two
outgoing edges: once both players have revealed their cards,
the winner, by putting the top cards to the bottom of his hand,
defines such two transitions with probabilities $p_1^i$ and
$p_2^i$, where $i$ is either $L$ or $R$, depending on who
wins this particular game. We call a vertex {\em attaining} if
it has a final state as one of its descendants and {\em
wandering} otherwise.  It is obvious that a descendant of a
wandering vertex is again wandering, and a predecessor of an
attaining one is again attaining. A graph is called {\em
absorbing} if all the vertices are attaining. That is, the
graph of our game is absorbing if and only if for every state
(each division on the pack into two hands) it is possible to
finish playing the game in a finite number of moves. The
difference between absorbing and non-absorbing graphs can be
seen in Fig.~\ref{fig:graph}.
Both graphs {\it i)} and {\it ii)} have the same final state $a$. Observe that the graph {\it i)} is not absorbing. It can be easily seen that it is impossible to reach the state $a$ from any of the vertices $c$, $d$ and $e$, while the graph {\it ii)}, which differs from {\it i)} only by one additional edge that leads from $e$ to $a$, is absorbing, as one can get into $a$ from any other vertex.

\section{Proof of the main result}

\begin{figure}[h]
\centering
\includegraphics[height=90pt, keepaspectratio]{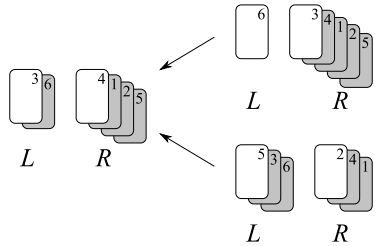}
 \caption{If the players have at
least two cards each, then this vertex has exactly two
predecessors.}\label{A2}
\end{figure}

We have already established that each state (except for final) has exactly two outgoing edges. Now assume that each of the players has at least two cards in their hands. There are exactly two incoming edges corresponding to the possible preceding  plays in which either of the players won (see Fig.~\ref{A2} for an illustration for a game with $6$ cards).

\begin{figure}[h]
\centering
\includegraphics[height=60pt, keepaspectratio]{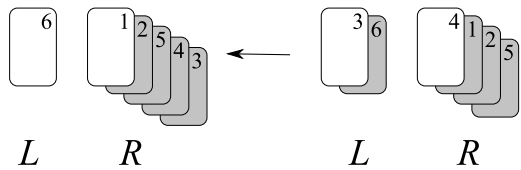}
\caption{If one of the players has only one card left,
this vertex has a unique direct predecessor.}\label{A1}
\end{figure}

If one of the players has only one card left, he could not have won the preceding play, and hence there is only one possibility for the winner. This is illustrated on Fig.~\ref{A1} for a game with $6$ cards.

Using this crucial observation, we are going to show that {\it the
game graph has no wandering vertices} (i.e. it is possible to
finish the game starting from any state). We show this by assuming
that there is at least one wandering vertex, and then establishing
that this leads to a contradiction.

By $\mathcal W$ denote the set of all wandering vertices in our
graph. Let $m$ be the number of vertices in $\mathcal W$. Every
vertice in $\mathcal W$ has exactly 2 descendants, so the
number of edges leaving from $\mathcal W$  is $2m$. Observe
that each edge going out from a wandering vertex goes into a
wandering vertex again, otherwise we would get a contradiction
with the definition of a wandering vertex.

Therefore, the total number of different edges that lead into
$\mathcal W$ is more or equal to  $2m$, since it includes all
the edges that leave $\mathcal W$ plus the edges from the
outside. Taking into account that each non-terminal vertex has
 either one or two incoming edges, we conclude that the amount of
the incoming edges could not be more than $2m$ so it is $2m$
exactly
 and every vertex from $\mathcal W$ has got exactly 2 incoming vertices. So we
immediately get two results:  If wandering
vertices exist, then they comprise an {\em isolated subgraph} (i.e. it is impossible to get into there from any vertex outside
of this subgraph); second, each wandering vertex has got
exactly two predecessors, and hence each wandering vertex
corresponds to the state in which each player has at least
two cards.

Now pick up any wandering vertex and conduct the following `back-tracking' procedure. If the vertex has two direct predecessors, consider the one in which the left player has less cards than in the current state (that is, the left player won the preceding play). Observe that this is always possible (see Fig.~\ref{A2}), as we can always backtrack through the play in which the left player was the winner. If we continue going back in this manner, we will finally reach the state in which the left player has  only one card left, and hence has  only one direct predecessor; this can not correspond to a wandering vertex. This means that a wandering vertex can be reached from a non-wandering, which contradicts our earlier finding (that wandering vertices constitute an isolated subgraph). This can only mean that the graph does not have any wandering vertices at all, and hence for each state there is a path to a final one.

Now we can use a well-established fact (see \cite[Chap.3]{fmc})
that if a graph corresponding to a finite Markov chain is
absorbing, the mathematical expectation of the number of moves
needed to reach the final state is finite, which is exactly what
we expected to prove in the first place.

 In particular, we prove that in the game with one
suit the player having the highest card wins in a finite time. Zdes sleduet otmetit, chto eto utverzhdenie ne verno v Russkoi versii igri, v kotoroi edinstvennim otlichiem yavlyaetsya to, chto starhaya karta proigrivaet samoi mladshei karte. Konechno, i dlya etoi versii veren nash, resultat o konechnosti igri. !!!!!!!!!

We will not discuss the proof for the classical game of war in
detail because the basic mathematical ideas are already
contained in the proof of the model game (For details
see \cite[p.8]{lr}) We only note that the main idea in studying
the real card game is the following obvious statement. {\em If
a subgraph of an oriented graph, which consists of all the
vertices of the original graph, and might not include some of
the edges, does not have any wandering vertices, then the
original graph does not have them either}.

We would like to point out that this card game was studied by
other authors, however, they focused on other aspects of the
game. In particular, Jacob Haqq-Misra \cite{art2} uses
numerical simulation (employing Monte-Carlo method) to find out
how the advantage in the initial distribution of cards
influences the outcome of the standard game; Ben-Naim and
Krapivsky \cite{art1} discuss a stochastic model of the game.
Very recently Michael Z. Spivey \cite{msp} has shown the
possibility of cycles in rather general assumptions about the
rules and the number of cards.

\begin{figure}[h]
\centering
\includegraphics[height=60pt, keepaspectratio]{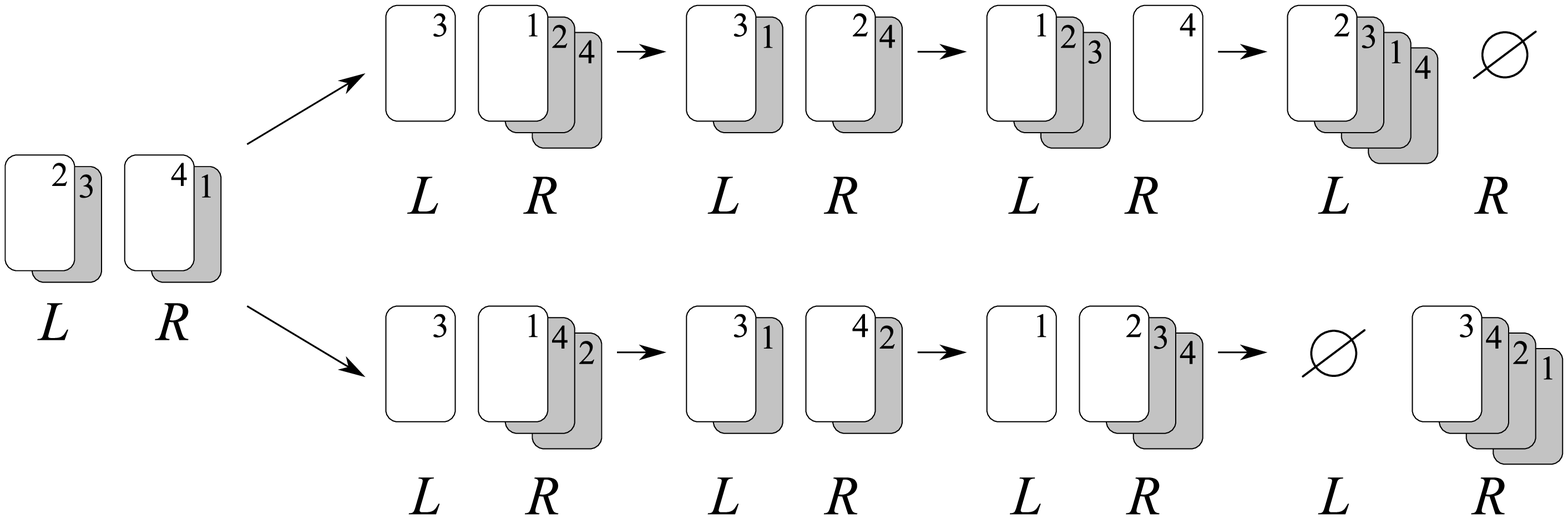}
\caption{A game with two possible outcomes: $n = 4$, the lowest valued card
in the pack wins over the highest valued one.}\label{F7}
\end{figure}

It is a curious observation that when the rules of the game
allow the lowest valued card in the pack to win over the
highest one, it is possible that with the same initial setting
either of the players can win. An example of such situation is
provided on Fig. \ref{F7}.

Finally let us note that for every $n$ the probability to keep
playing after $n$ steps is always positive. It follows from the
fact that the graph of the game has cycles.\footnote{We suppose
that initial state is arbitrary}  To show this, it is enough to
pick any non-final state and transition to any of the
ancestors. The existence of a cycle follows from the finiteness
of the graph.

\paragraph{Acknowledgments.}
The authors are grateful
to the first author's son Grigory Lakshtanov (5 years old) for tireless practical experimentation searching for cycles with his father, and to the second author's partner Ian Robson for
his help with correcting the authors' English.

\bigskip

\noindent\textbf{Evgeny Lakshtanov} received his Ph.D. from the Moscow State University 
in December 2004. He currently does his posdoctoral research in the University of Aveiro, Portugal. On schitaet samim vazhnim resultatom svoei raboti sozdanie of simplified versions of  serious mathematical models to make them understandable by children, preserving both aesthetic and intellectual value. The latter is in particularly measured by whether a given simplification allows setting a sufficient list of problems feasible for school students.

\noindent\textit{ Center for Research and Development in
Mathematics and Applications (CIDMA), Department of
Mathematics, Aveiro University,  Portugal.  \\
lakshtanov@ua.pt}

\bigskip
\noindent\textbf{Vera Roshchina} received her... When she is not doing
mathematics, she enjoys singing, bicycling, and playing volleyball

\noindent\textit{Research Center in Mathematics and
Applications (CIMA), University of Evora, Portugal. \\ }

\end{document}